\documentclass{amsart}
\usepackage{amsmath}
\usepackage{amssymb}
\usepackage{amscd}
\usepackage[all]{xy}
\usepackage{color}
\usepackage{fancyhdr}

\theoremstyle{plain}
\newtheorem{theorem}{Theorem}

\newtheorem{definition}[theorem]{Definition}

\theoremstyle{remark}
\newtheorem{remark}[theorem]{Remark}
\newtheorem{example}[theorem]{Example}

\def\bc{\begin{center}}
\def\ec{\end{center}}

\def\Spec{{\rm Spec}}

\def\ch{{\rm char}}

\begin{document}
\title[almost discrete valuation domains] {almost discrete valuation domains}

\author [Anderson] {D. D. Anderson$^\ast$}
\address{(Anderson) Department of Mathematics, The University of Iowa, Iowa City, IA 52242-1419, USA}
\email{dan-anderson@uiowa.edu}

\author [Xing] {Shiqi Xing}
\address{(Xing) College of Applied Mathematics, Chengdu University of Information Technology, Chengdu, Sichuan 610225, China}
\email{sqxing@yeah.net}

\author [Zafrullah] {Muhammad Zafrullah}
\address{(Zafrullah) Department of Mathematics, Idaho State University, Pocatello, ID 83209-8085, USA}
\email{mzafrullah@usa.net}


\thanks{$2010$ Mathematics Subject Classification: 13A15, 13A18, 13F05, 13F25, 13G05, 13J05, 13J10}

\thanks{$\ast$ Corresponding Author}

\date{\today}

\begin{abstract} Let $D$ be an integral domain. Then $D$ is an almost valuation (AV-)domain if for $a, b\in D\setminus \{0\}$ there exists a natural number $n$ with $a^{n}\mid b^{n}$ or $b^{n}\mid a^{n}$. AV-domains are closely related to valuation domains, for example, $D$ is an AV-domain if and only if the integral closure $\overline{D}$ is a valuation domain and $D\subseteq \overline{D}$ is a root extension. In this note we explore various generalizations of DVRs (which we might call almost DVRs) such as Noetherian AV-domains, AV-domains with $\overline{D}$ a DVR, and quasilocal and local API-domains (i.e., for $\{a_{\alpha}\}_{\alpha\in \Lambda}\subseteq D$, there exists an $n$ with $(\{a_{\alpha}^{n}\}_{\alpha\in \Lambda})$ principal). The structure of complete local AV-domains and API-domains is determined.
\end{abstract}

\maketitle

Let $D$ be a commutative integral domain. Recall from \cite[Definition 5.5]{AZ02} that a domain $D$ is called an \textit{almost valuation domain} (\textit{AV-domain}) if for $a, b\in D^{\ast}:=D\setminus\{0\}$, there exists a natural number $n=n(a,b)$ with $a^{n}\mid b^{n}$ or $b^{n}\mid a^{n}$. As expected, AV-domains are closely related to valuation domains. For example (see Theorem \ref{02}), a domain $D$ is an AV-domain if and only if $\overline{D}$, the integral closure of $D$, is a valuation domain and $D\subseteq \overline{D}$ is a root extension (i.e., for each $x\in \overline{D}$, there exists a natural number $n=n(x)$ with $x^{n}\in D$). The purpose of this article is to study AV-domains which in some way generalize the notion of a DVR. (Here by a DVR (discrete valuation ring) we mean a rank-one discrete valuation domain, or equivalently, a local PID which is not a field.) As we shall see there are plenty of different possible ways to define an ``almost DVR". But first we give a little background.

In his study of almost factoriality, Zafrullah \cite{Z} introduced the notion of an almost GCD domain (AGCD-domain). Recall that a domain $D$ is called an \textit{AGCD-domain} if for $a,b\in D^{\ast}$, there exists a natural number $n=n(a,b)$ with $a^{n}D\bigcap b^{n}D$ (or equivalently, $(a^{n},b^{n})_{v}$) principal. It was shown in \cite[Theorems 3.1 and 3.4]{Z} that (1) if a domain $D$ is an AGCD-domain, then $D\subseteq \overline{D}$ is a root extension and $\overline{D}$ is an AGCD-domain and (2) an integrally closed domain is an AGCD-domain if and only if it is a P$v$MD with torison $t$-class group \cite[Corollary 3.5 and Theorem 3.9]{Z}.

Anderson and Zafrullah \cite{AZ02} introduced the notion of an almost B\'{e}zout domain (AB-domain) and these domains were further studied in \cite{AKL}. Recall that $D$ is called an \textit{AB-domain} if for $a,b\in D^{\ast}$, there exists a natural number $n=n(a,b)$ with $(a^{n}, b^{n})$ principal. It was shown that (1) a domain $D$ is an AB-domain if and only if $D\subseteq \overline{D}$ is a root extension and $\overline{D}$ is an AB-domain \cite[Theorem 4.6]{AZ02} and (2) an integrally closed domain is an AB-domain if and only if it is a Pr\"{u}fer domain with torsion class group \cite[Theorem 4.7]{AZ02}. Also a domain $D$ is called an \textit{almost principal ideal domain} (\textit{API-domain}) if for any nonempty subset $\{d_{\alpha}\}_{\alpha\in \Lambda}$ of $D^{\ast}$, there exists a natural number $n=n(\{d_{\alpha}\}_{\alpha\in \Lambda})$ with $(\{d_{\alpha}^{n}\}_{\alpha\in \Lambda})$ principal. Now an integrally closed domain is an API-domain if and only if it is a Dedekind domain with torsion class group \cite[Theorem 4.12]{AZ02}. However, unlike the AV-domain and AB-domain cases, we can have $D\subseteq \overline{D}$ a root extension and $\overline{D}$ is an API-domain (even a DVR) without $D$ being an API-domain \cite[Example 4.14]{AZ02}. Of course an API-domain is an AB-domain and hence $D\subseteq \overline{D}$ is a root extension with $\overline{D}$ being an AB-domain (equivalently, a Pr\"{u}fer domain with torsion class group), but it is not known whether $\overline{D}$ must actually be a Dedekind domain. In the case where $D\subseteq \overline{D}$ is a bounded root extension (i.e., there exists a natural number $n$ with $x^{n}\in D$ for each $x\in \overline{D}$), $D$ is an API-domain if and only if $\overline{D}$ is an API-domain \cite[Theorem 4.11]{AZ02}. However, if $D\subseteq \overline {D}$ is a root extension and $\overline{D}$ is an AGCD-domain (even a UFD), $D$ need not be an AGCD-domain (see \cite[Example 3.1]{AKL}).

We next give a brief review of root extensions. An extension $R\subseteq S$ of commutative rings is a \textit{root extension} if for each $s\in S$, there exists a natural number $n=n(s)$ with $s^{n}\in R$ and is a \textit{bounded root extension} if there exists a natural number $n$ with $s^{n}\in R$ for each $s\in S$. If $R\subseteq S$ is a root extension, then the map $\Spec (S)\rightarrow \Spec(R)$ given by $Q\rightarrow Q\bigcap R$ is a homeomorphism with inverse $P\rightarrow \{s\in S\mid s^{n}\in P$ for some natural number $n\}$ \cite[Theorem 2.1]{AZ02}. Hence $\dim R=\dim S$ (this also follows since $R\subseteq S$ is an integral extension) and $R$ is quasilocal if and only if $S$ is quasilocal. Note that if $R\subseteq S$ is a (bounded) root extension, then $T(R)\subseteq T(S)$ is also a (bounded) root extension, where $T(R)$ (resp., $T(S)$) is the total quotient ring of $R$ (resp., $T(S)$). The following theorem characterizes when a field extension is a (bounded) root extension.

\begin{theorem}\label{01}
Let $K\subseteq L$ be a field extension.
\begin{itemize}
\item[(1)] (\cite{N}) $K\subseteq L$ is a root extension if and only if $L/K$ is purely inseparable (which in the characteristic 0 case means $K=L$) or $L$ is algebraic over a finite field.
\item[(2)] (\cite[Corollary 2.2]{AKL}) $K\subseteq L$ is a bounded root extension if and only if $L$ is purely inseparable over $K$ of bounded exponent (i.e, $K=L$ or $L^{p^{n}}\subseteq K$ for some natural number $n$ where $\ch K=p$) or $L$ is finite.
\item[(3)]  (\cite[Corollary 2.3]{AKL}) $K\subseteq L$ is a root extension with the property that for each intermediate field $F$ ($K\subseteq F\subseteq L$) with $[F:K]$ finite, $K\subseteq F$ is a bounded root extension if and only if $L/K$ is purely inseparable or $K$ is a finite field and $L$ is algebraic over $K$.
\end{itemize}
\end{theorem}

We next give some results about AV-domains, some known and some new.
\begin{theorem}\label{02}
\begin{itemize}
\item[(1)] For an integral domain $D$ with quotient field $K$, the following conditions are equivalent.
\begin{itemize}
\item[(a)] $D$ is an AV-domain.
\item[(b)] $\overline{D}$ is a valuation domain and $D\subseteq \overline{D}$ is a root extension.
\item[(c)] $D$ is a $t$-local AGCD-domain (i.e., $D$ has a unique maximal $t$-ideal).
\item[(d)] $D$ is a quasilocal AB-domain.
\item[(e)] If $0\neq x\in K^{\ast}$, there exists a natural number $n=n(x)$ with $x^{n}\in D$ or $x^{-n}\in D$.
\end{itemize}
\item[(2)] A root closed domain $D$ (and hence an integrally closed domain $D$) is an AV-domain if and only if $D$ is a valuation domain.
\item[(3)] An overring of an AV-domain is an AV-domain.
\item[(4)] Let $R\subseteq S$ be a root extension of integral domains. Then $R$ is an AV-domain if and only if $S$ is an AV-domain.
\item[(5)] Let $D$ be an AV-domain and $P$ a prime ideal of $D$. Then $D/P$ is an AV-domain.
\item[(6)] Let $D$ is an AV-domain with the quotient field $K$ and let $L$ be a subfield of $K$. Then $D\bigcap L$ is an AV-domain. Let $M$ be the quotient field of $D\bigcap L$. Then $\overline{D}\bigcap M$ is an valuation domain with quotient field $M$ and $\overline{D\bigcap L}=\overline{D\bigcap M}=\overline{D}\bigcap M$.
\end{itemize}
\end{theorem}
\begin{proof}
$(1)$ \ \ The equivalence $(a)$-$(d)$ is \cite[Theorem 5.6]{AZ02}, while the equivalence $(a)\Leftrightarrow(e)$ is clear.

$(2)$ \ \ This immediately follows from (1).

$(3)$ \ \ This immediately follows from the equivalence $1(a)$ and $1(e)$.

$(4)$\ \ $(\Rightarrow)$ Let $a, b\in S^{\ast}$, so there exists a natural numbers $n$ with $a^{n}, b^{n}\in R$. So there exists a natural number $m$ with $a^{nm}\mid b^{nm}$ or $b^{nm}\mid a^{nm}$ in $R$. Hence $a^{nm}\mid b^{nm}$ or $b^{nm}\mid a^{nm}$ in $S$.

$(\Leftarrow)$ Let $a, b\in R^{\ast}$. So there exists a natural number $n$ with $a^{n}\mid b^{n}$ or $b^{n}\mid a^{n}$ in $S$. Without loss of generality, suppose that $a^{n}\mid b^{n}$ in $S$, so $sa^{n}=b^{n}$ for some $s\in S$. Choose $m$ with $s^{m}\in D$. Then $s^{m}a^{mn}=b^{mn}$. Hence $a^{mn}\mid b^{mn}$ in $R$.

$(5)$ \ \ Suppose that $D$ is an AV-domain. Let $x, y\in D/P$, so $x=a+P$ and $y=b+P$ for some $a, b\in D$.  Without loss of generality, we can assume that $(a^{n})\subseteq (b^{n})$ for some natural number $n$. Then $(x^{n})\subseteq (y^{n})$. So $D/P$ is an AV-domain.

$(6)$\ \ Now $\overline{D}$ is a valuation domain with quotient field $K$. So $\overline{D}\bigcap L$ is a valuation domain with field $L$. Let $M$ be the quotient field of $D\bigcap L$ (unlike the valuation domain case, $M$ may be a proper subfield of $L$). Then $D\bigcap L=D\bigcap M$ and $\overline{D}\bigcap M$ is a valuation overring of $D\bigcap L$ with quotient field $M$. Since $D\subseteq \overline{D}$ is a root extension, $D\bigcap M\subseteq \overline{D}\bigcap M$ is a root extension. Then $\overline{D\bigcap M}=\overline{D}\bigcap M$ and $D\bigcap M$ is an AV-domain.
\end{proof}
\remark \label{03}
\begin{itemize}
\item[(1)]\ \ Let $D$ be a subring of a field $K$. If for each $x\in K^{\ast}$, either $x\in D$ or $x^{-1}\in D$, then $D$ is a valuation domain with quotient field $K$. Suppose that we just have for each $x\in K^{\ast}$, that $x^{n}\in D$ or $x^{-n}\in D$ for some natural number $n=n(x)$. Then while $D$ is an AV-domain, $K$ need not be the quotient field of $D$. For example, let $M\subsetneq N$ be a purely
inseparable field extension of bounded exponent. Take $D=M[[X]]$ and $K=N((X))$, the quotient field of $N[[X]]$. For $f\in K^{\ast}$, there is a natural number $n=n(f)$ with $f^{n}\in M[[X]]$ or $f^{-n}\in M[[X]]$, but the quotient field of $D$ is $M((X))\subsetneq K$.

\item [(2)]\ \ If $V$ is a valuation domain with quotient field $K$ and $L$ is a subfield of $K$, then $V\bigcap L$ is a valuation domain with quotient field $L$. However, if $V$ is an AV-domain, $V\bigcap L$ need not have quotient field $L$. Let $F\subsetneq G$ be a purely
inseparable field extension. Let $R=F+XG[[X]]$, so $\overline{R}=G[[X]]$ a DVR. Now $R\subseteq \overline{R}$ is a root extension, so $R$ is an AV-domain with quotient field $G((X))$. But $R\bigcap F=F\subsetneq G=\overline{R}\bigcap G$.\\
\end{itemize}

We next collect some results about API-domains.\\

\begin{theorem}\label{04}
The following statements hold for API-domains.
\begin{itemize}
\item[(1)] A quasilocal API-domain is an AV-domain.
\item[(2)] A domain $D$ is a quasilocal API-domain if and only if for each nonempty subset $\{a_{\alpha}\}_{\alpha\in\Lambda}$ of $D^{\ast}$, there exists a natural number $n=n(\{a_{\alpha}\}_{\alpha\in\Lambda})$ and an $\alpha_{0}\in \Lambda$, with $a_{\alpha_{0}}^{n}\mid a_{\alpha}^{n}$ for each $\alpha\in \Lambda$.
\item[(3)] A root closed domain $D$ (and hence an integrally closed domain $D$) is an API-domain if and only if $D$ is a Dedekind domain with torsion class group. Thus a root closed quasilocal domain $D$ (and hence an integrally closed quasilocal domain $D$) is an API-domain if and only if $D$ is a DVR.
\item[(4)] Let $R\subseteq S$ be a bounded root extension of domains. Then $R$ is a (quasilocal) API-domain if and only if $S$ is a (quasilocal) API-domain.
\item[(5)] Let $D$ be an API-domain with $[D:\overline{D}]\neq 0$. Then $D\subseteq \overline{D}$ is a bounded root extension and hence $\overline{D}$ is an API-domain. Thus if $D$ is a quasilocal API-domain with $[D:\overline{D}]\neq 0$, $\overline{D}$ is a DVR.
\item[(6)] Let $D$ be an API-domain and $S$ a multiplicatively closed subset of $D$. Then $D_{S}$ is an API-domain. If further $D$ is quasilocal, so is $D_{S}$.
\item[(7)] $D$ is an API-domain if and only if $D_{M}$ is an API-domain for each maximal ideal $M$ of $D$ and $D$ has torsion $t$-class group $Cl_{t}(D)$.
\item[(8)]If $D$ is a (quasilocal) API-domain and $P$ is a prime ideal of $D$, then $D/P$ is a (quasilocal) API-domain.
\item[(9)] Let $D$ be an API-domain containing a field of characteristic 0. Then $D$ is a Dedekind domain with torsion class group. Hence if $(D, M)$ is a quasilocal API-domain with $\mbox{char} D/M=0$, $D$ is a DVR.
\item[(10)] Let $D$ be a quasilocal API-domain with the quotient field $K$. Let $L$ be a subfield of $K$. Then $D\bigcap L$ is a quasilocal API-domain.
\end{itemize}
\end{theorem}
\begin{proof}
$(1)$\ \ This follows from Theorem \ref{02} (1) or from (2).

$(2)$\ \ $(\Leftarrow)$ Certainly $D$ is an AV-domain, so $D$ is quasilocal. Note that $(\{a^{n}_{\alpha}\}_{\alpha\in\Lambda})=(a^{n}_{\alpha_{0}})$, so $D$ is an API-domain.

$(\Rightarrow)$ Since $D$ is an API-domain, $(\{a_{\alpha}^{n}\}_{\alpha\in\Lambda})$ is principal for some natural number $n$. Since $D$ is quasilocal, a principal ideal is completely join-irreducible. Thus $(\{a_{\alpha}^{n}\}_{\alpha\in\Lambda})=(a_{\alpha_{0}}^{n})$ for some $\alpha_{0}\in \Lambda$.

$(3)$\ \ \cite[Theorem 4.12]{AZ02}.

$(4)$\ \ \cite[Theorem 4.11]{AZ02} gives that $R$ is an API-domain if and only if $S$ is an API-domain. Since $R\subseteq S$ is a root extension, $R$ is quasilocal if and only if $S$ is quasilocal.

$(5)$\ \ Suppose that $D$ is an API-domain with $[D:\overline{D}]\neq 0$. Let $0\neq m\in [D:\overline{D}]$. Then $\{md\mid d\in \overline{D}^{\ast}\} \subseteq D^{\ast}$, so there exists a natural number $n$ with $(\{m^{n}d^{n}\mid d\in \overline{D}^{\ast}\})$ a principal ideal of $D$. Hence there exist $d_{1},\dots, d_{s}\in \overline{D}^{\ast}$ with $(\{m^{n}d^{n}\mid d\in \overline{D}^{\ast}\})= Dm^{n}d_{1}^{n}+\dots+Dm^{n}d_{s}^{n}$. Since $D\subseteq \overline{D}$ is root extension, there exists a natural number $k$ with $d_{i}^{nk}\in D$ for $i=1,\dots, s$. Since $(\{m^{n}d^{n}\mid d\in \overline{D}^{\ast}\})= m^{n}d_{1}^{n}D+\dots+m^{n}d_{s}^{n}D$ are principal, $(\{m^{n}d^{n}\mid d\in \overline{D}^{\ast}\})^{k}=(\{m^{nk}d^{nk}\mid d\in \overline{D}^{\ast}\})^{k}=(m^{n}d_{1}^{n}D+\dots+m^{n}d_{s}^{n}D)^{k}=m^{nk}d_{1}^{nk}D+\dots+m^{nk}d_{s}^{nk}D$. Hence $d^{nk}m^{nk}\in m^{nk}d_{1}^{nk}D+\dots+m^{nk}d_{s}^{nk}D$, so $d^{nk}\in Dd_{1}^{nk}+\dots+Dd_{s}^{nk}\subseteq D$. Thus $D\subseteq \overline {D}$ is a bounded root extension. Suppose that $D$ is quasilocal. Then $\overline{D}$ is an integrally closed quasilocal API-domain. So $\overline{D}$ is a quasilocal Dedekind domain and hence $\overline{D}$ is a DVR.

$(6)$ \ \ Let $\{d_{\alpha}/s_{\alpha}\}_{\alpha\in\Lambda}\in D^{\ast}_{S}$ where $d_{\alpha}\in D^{\ast}$ and $s_{\alpha}\in S$. Since $D$ is an API-domain, there is a natural number $n$ with $(\{d_{\alpha}^{n}\}_{\alpha\in\Lambda})$ a principal ideal of $D$. Then $D_{S}(\{(d_{\alpha}/s_{\alpha})^{n}\}_{\alpha\in \Lambda})=(\{d_{\alpha}^{n}\}_{\alpha\in\Lambda})_{S}$ is a principal ideal of $D_{S}$. So $D_{S}$ is an API-domain. The last statement follows from the fact that a quasilocal API-domain is an AV-domain and hence $\Spec (D)$ is totally ordered.

$(7)$\ \ The proof that $D$ is an AP-domain if and only if each $D_{M}$ is an AV-domain \cite[Theorem 5.8]{AZ02} carries over mutatis mutandis to give that a domain $D$ is an AD-domain (i.e, for $\{x_{\alpha}\}_{\alpha\in \Lambda}\subseteq D^{\ast}$ there exists a natural number $n$ with $(\{x_{\alpha}^{n}\}_{\alpha\in \Lambda})$ invertible) if and only if $D_{M}$ is an API-domain for each maximal ideal $M$ of $D$. Now an API-domain being an AGCD-domain has a torsion $t$-class group \cite[Theorem 3.4]{AZ02}. But an AD-domain with torsion $t$-class group is certainly an API-domain \cite[Theorem 4.11]{AZ02}.

$(8)$\ \ The proof is similar to the proof of Theorem \ref{02} (5).

$(9)$\ \ Let $I$ be a nonzero ideal of $D$. Then there exists a natural number $n$ with $I_{n}:=(\{i^{n}\mid i\in I\})$ principal. Since $D$ contains a field of characteristic 0, $I_{n}=I^{n}$ \cite[Theorem 6.12]{AZ02}. So $I^{n}$ is principal and hence $I$ is invertible. So $D$ is a Dedekind domain with torsion class group. Note that a quasilocal domain $(D,M)$ contains a field of characteristic 0 if and only if $\ch D/M=0$. The last statement is now immediate.

$(10)$\ \ Let $\{r_{\alpha}\}_{\alpha\in \Lambda}\subseteq D^{\ast}\bigcap L$.  Since $D$ is a quasilocal API-domain, by (2) there exists a natural number $n$ and $\alpha_{0}\in \Lambda$ with $r_{\alpha_{0}}^{n}\mid r_{\alpha}^{n}$ in $D$ for each $\alpha\in \Lambda$. Now $r_{\alpha_{0}}^{n}d_{\alpha}=r_{\alpha}^{n}$ for some $d_{\alpha}\in D$. Hence
$d_{\alpha}=r_{\alpha}^{n}/r_{\alpha_{0}}^{n}\in D\bigcap L$. So $r_{\alpha_{0}}^{n}\mid r_{\alpha}^{n}$ in $D\bigcap L$.
\end{proof}
We next note that the converse of Theorem \ref{04} (5) is not true. $D$ can be a local API-domain with $D\subseteq \overline{D}$ a bounded root extension, but $[D:\overline{D}]=0$ (equivalently, $\overline{D}$ is not a finitely generated $\overline{D}$-module).

\begin{example}\label{0004}
($D$ a local API-domain with $D\subseteq \overline{D}$ a bounded root extension $\nRightarrow$ $[D:\overline{D}]\neq 0$.) This is Nagata's example of a
one-dimensional local domain $D$ with $\overline{D}$ not a finitely generated $D$-module \cite[E 3.2, page 206]{N01}.) \rm{Let $K=\mathbb{Z}_{p}(\{Y_{i}\}_{i=1}^{\infty})$ where $\{Y_{i}\}_{i=1}^{\infty}$ is a set of indeterminates over $\mathbb{Z}_{p}$, so $K^{p}=\mathbb{Z}_{p}(\{Y_{i}^{p}\}_{i=1}^{\infty})$. Let $R=K^{p}[[X]][K]$ where $X$ is a power series indeterminate over $K^{p}$ and $c=\sum_{i=1}^{\infty}Y_{i}X^{i}$. Put $D=R[c]$. So $\overline{D}$ is a DVR and $\overline {D}$ is not a finitely generated $D$-module. Note that for $f\in D$ (resp., $f\in T(D)$), $f^{p}\in K^{p}[[X^{p}]]$ (resp., $f^{p}\in K^{p}((X^{p}))$). Let $z\in \overline{D}$, so $z^{n}+r_{1}z^{n-1}+\dots+r_{n}=0$ where $r_{i}\in D$. Then $(z^{p})^{n}+r_{1}^{p}(z^{p})^{n-1}+\dots+r_{n}^{p}=0$ where $r_{i}^{p}\in D$. So $z^{p}\in K^{p}((X^{p}))$ is integral over $Z_{p}[[X^{p}]]$ and hence $z^{p}\in \mathbb{Z}_{p}[[X^{p}]]\subseteq D$. So $D\subseteq \overline{D}$ is a bounded root extension.
}
\end{example}

Let $R$ be a commutative ring. The lattice $L(R)$ of ideals of $R$ is a complete multiplicative lattice. Let $R$ and $S$ be two quasilocal rings and suppose that $\Theta: L(R)\rightarrow L(S)$ is a complete multiplicative lattice isomorphism. Recall that in a quasilocal ring an ideal is principal if and only if it is completely join-irreducible. Thus an ideal $J$ of $R$ is principal if and only if $\Theta (J)$ is principal. It easily follows that $R$ is an AV-domain (resp., API-domain) if and only if $S$ is an AV-domain (resp., API-domain). (The hypothesis that $R$ and $S$ are quasilocal is essential. Indeed, for any Pr\"{u}fer domain (resp., Dedekind domain) $D$, $D(X)$ is a Bezout domain (resp, PID) and the map $\Theta: L(D)\rightarrow L(D(X))$ given by $\Theta(J)=D(X)J$ is a complete multiplicative lattice isomorphism, c.f. \cite[Theorem 8]{A}.)

Suppose that $(D,M)$ is a one-dimensional local domain. Then the map $\Theta: L(D)\rightarrow L(\widehat{D})$ given by $\Theta(J)=\widehat{D}J$, where $\widehat{D}$ is the $M$-adic completion of $D$, is a complete multiplicative lattice isomorphism if and only if $\widehat{D}$ is analytically irreducible (i.e., $\widehat{D}$ is an integral domain). The implication ($\Rightarrow$) is clear. Conversely, suppose that $\widehat{D}$ is an integral domain. Then any ideal of $D$ (resp., $\widehat{D}$) other than $0$ and $D$ (resp., $\widehat{D}$) is $M$-primary (resp., $\widehat{M}$--primary). Since the map $Q\rightarrow \widehat{D}Q$ (with inverse the contraction map) is an order preserving bijection between the set of $M$-primary ideals of $D$ and the set of $\widehat{M}$-primary ideals of $\widehat{D}$, $\Theta$ is a a complete multiplicative lattice isomorphism. Thus for $D$ analytically irreducible, $D$ is an AV-domain (resp., API-domain) if and only if $\widehat{D}$ is an AV-domain (resp., API-domain). Recall that $D$ is analytically irreducible if and only if $\overline{D}$ is finitely generated $D$-module (equivalently, $[D:\overline{D}]\neq 0$) and $\overline{D}$ is a DVR \cite{N01}. So Example \ref{0004} shows that a local API-domain $D$ need not be analytically irreducible and hence $\widehat{D}$ need not be an API-domain. Thus we have the following theorem.
\begin{theorem}\label{004}
Let $(D, M)$ be a local domain. Consider the following statements.
\begin{itemize}
\item[(1)] $D$ is an AV-domain (resp., API-domain).
\item[(2)] $D\subseteq \overline{D}$ is a root extension (resp., bounded root extension) and $\overline{D}$ is a DVR (equivalently, $\dim D=1$).
\item[(3)] $\widehat{D}$ is an AV-domain (resp., API-domain).
\end{itemize}
Then $(3)\Rightarrow (2)\Rightarrow (1)$ and if $\overline {D}$ is a finitely generated $D$-module, $(1)\Rightarrow (3)$.
\end{theorem}
\begin{proof}
Since $(1), (2)$ and $(3)$ give that $\dim D=1$ and for a root extension $D\subseteq \overline{D}$, $\overline{D}$ is a DVR if and only if $\dim D=1$, we may assume that $\dim D=1$ and replace $(2)$ by $(2')$: $D\subseteq \overline{D}$ is a root extension (resp., bounded root extension).

$(3)\Rightarrow (2)$\ \ Suppose that $\widehat{D}$ is an AV-domain (resp., API-domain). So $\widehat{D}$ is an integral domain and hence from the previous paragraph, $\widehat{D}$ an AV-domain (resp., API-domain) gives that $D$ is an AV-domain (resp., API-domain). Thus if $D$ is an AV-domain, $D\subseteq \overline {D}$ is a root extension. Suppose that $\widehat{D}$ is an API-domain. Then $\widehat{D}$ being a domain gives $[D:\overline{D}]\neq 0$, and hence by Theorem \ref{04} (5), $D\subseteq \overline{D}$ is a bounded root extension.

$(2)\Rightarrow (1)$ \ \ Suppose that $D\subseteq \overline{D}$ is a root extension. Then $\dim D=1$ gives that $\overline{D}$ is a DVR and hence $D$ is an AV-domain. Suppose that $D\subseteq \overline{D}$ is a bounded root extension. Then $D$ is an API-domain by Theorem \ref{04} (4).

Suppose that $\overline{D}$ is a finitely generated $D$-module. If $D$ is an AV-domain (or an API-domain), then $\overline{D}$ is a DVR. Hence $D$ is analytically irreducible. By the remarks of paragraph preceding Theorem \ref{004}, we have that $D$ an AV-domain (resp., API-domain) gives that $\widehat{D}$ is an AV-domain (resp., API-domain).
\end{proof}

The next two theorems give a complete characterization of complete local AV-domains and API-domains.

\begin{theorem}\label{005}
(Equicharacteristic case)
\begin{itemize}
\item[(1)] Let $(D, M)$ be a complete local AV-domain (resp., API-domain) with $\ch$ $D/M(=\ch D)=0$. Then $D$ is a DVR and hence $D\cong K[[X]]$ where $K\cong D/M$ is a coefficient field for $D$.
\item[(2)] Let $k\subseteq K$ be a root extension of fields (resp., bounded root extension) where $\ch k=p>0$. Let $m$ be a natural number. Suppose that $D$ is a domain with $k+K[[X]]X^{m}\subseteq D\subseteq K[[X]]$. Then $D$ is a quasilocal AV-domain (resp., quasilocal API-domain) with $\overline{D}=K[[X]]$. More precisely, $(D,M)$ is an API-domain if and only if $D/M\subseteq K$ is a bounded root extension. Suppose that $[K:k]<\infty$. Then $D$ is a complete local AV-domain (resp., complete local API-domain).

\item[(3)] Suppose that $(D, M)$ is a complete local AV-domain (resp., complete local API-domain) with $\ch D=\ch D/M=p>0$. Then $(\overline{D}, \overline{M})$ is a complete local DVR, so $\overline{D}\cong K[[X]]$ where $K\cong \overline{D}/\overline{M}$ is a coefficient field for $\overline{D}$. Let $k\cong D/M$ be a coefficient field for $D$. So $D\subseteq \overline{D}$ and hence $k\subseteq K$ is a root extension (resp., bounded root extension) with $[K:k]<\infty$. Moreover, $k+K[[X]]X^{m}\subseteq D\subseteq K[[X]]$ where $0\neq [D:\overline{D}]=K[[X]]X^{m}$.
\end{itemize}
\end{theorem}
\begin{proof} $(1)$\ \ Suppose $(D, M)$ is a complete local AV-domain. Now $(\overline{D}, \overline{M})$ is a complete DVR and $D\subseteq \overline{D}$ is a root extension. Thus $D/M\subseteq \overline{D}/\overline{M}$ is a root extension. Since $\ch D/M=0$, $D/M=\overline{D}/\overline{M}$. Choose a coefficient field $K$ for $D$ (and hence for $\overline{D}$). So we can take $\overline{D}=K[[X]]$. Since $\overline{D}$ is a finitely generated $D$-module, there exists a natural number $n$ with $0\neq [D:\overline{D}]=K[[X]]X^{n}$. So $K+K[[X]]X^{n}\subseteq D\subseteq K[[X]]$ and $D\subseteq K[[X]]$ is a root extension. Suppose that $n>1$. Now there exists a natural number $m\geq 2$ with $(1+X^{n-1})^{m}\in D$. Now $(1+X^{n-1})^{m}-(1+mX^{n-1})=(^m_2)(X^{n-1})^{2}+\dots+(X^{n-1})^{m}\in K[[X]]X^{n}\subseteq D$, so $1+mX^{n-1}\in D$ and hence $X^{n-1}\in D$. Thus $K[[X]]X^{n-1}\subset D$, a contradiction. Thus $n=1$, so $D=K[[X]]$. Suppose that $(D,M)$ is a complete local API-domain. Then $D$ is also a complete local AV-domain and hence $D\cong K[[X]]$. However, in this case we can also just quote Theorem \ref{04} (9).

$(2)$\ \ Let $R=k+K[[X]]X^{m}$, so $R\subseteq D\subseteq K[[X]]$ with $\overline{R}=\overline{D}=K[[X]]$. Choose $l\geq 1$ with $p^{l}\geq m$. Let $f=\sum_{i=0}^{\infty}a_{i}X^{i}\in K[[X]]$, so $f^{p^{l}}=\sum_{i=1}^{\infty}a_{i}^{p^{l}}X^{p^{l}}$. Suppose that $a_{0}^{n}\in k$. Then $f^{p^{l}n}=(f^{p^{l}})^{n}=(a_{0}^{n})^{p^{l}}+n(a_{0}^{p^{l}})^{n-1}X^{p^{l}}+\dots \in R$. Hence if $k\subseteq K$ is a root extension (resp., bounded root extension), $R\subseteq K[[X]]$ is a root extension (resp., bounded root extension). Thus $D\subseteq K[[X]]$ is a root extension (resp., bounded root extension). So $D$ is an AV-domain (resp., quasilocal API-domain). Suppose that $D$ is an API-domain. Since $[D:\overline{D}]\neq 0$, Theorem \ref{04} (5) gives that $D\subseteq \overline{D}$ is a bounded root extension. Hence $D/M\subseteq K$ is a bounded root extension. Conversely, suppose that $D/M\subseteq K$ is a bounded root extension. Let $f\in \sum_{i=0}^{\infty}a_{i}X^{i}\in K[[X]]$. Choose $n$ with $x^{n}\in D/M$ for each $x\in K$. So there exists a $d\in D$ with $a_{0}^{n}-d\in K[[X]]X$. So $a_{0}^{np^{l}}-d^{p^{l}}\in D$. Thus $f^{p^{l}n}\in D$, so $D\subseteq \overline{D}$ is a bounded root extension. Thus $D$ is an API-domain by Theorem \ref{04} (4). Suppose that $[K: k]<\infty$. Then $K[[X]]$ is a finitely generated $R$-module and hence $R$ is Noetherian by the Eakin-Nagata Theorem. But $D$ is a finitely generated $R$-module, so $D$ is also Noetherian. Also, $R$ is a complete local ring and hence so is $D$ since $D$ is a finitely generated $R$-module.

$(3)$\ \ Suppose that $(D,M)$ is a complete local AV-domain (resp., complete local API-domain) with $\ch D=\ch D/M=p>0$. Now $(\overline{D}, \overline{M})$ is a complete DVR with $\ch \overline{D}=\ch \overline{D}/\overline{M}=p>0$, so $\overline{D}\cong K[[X]]$ where $K$ is a coefficient field for $\overline{D}$. Moreover, we can take a coefficient field $k\cong D/M$ for $D$ with $k\subseteq K$. Since $D$ is complete, $\overline{D}$ is a finitely generated $D$-module. So $0\neq [D:\overline {D}]=K[[X]]X^{m}$ for some natural number $m$. Since $D\subseteq \overline{D}$ is a root extension (resp., bounded root extension), $k\subseteq \overline{K}$ is a root extension (resp., bounded root extension). Note that $k+K[[X]]X^{m}\subseteq D\subseteq K[[X]]$. Since $\overline{D}$ is a finitely generated $D$-module, $\overline{D}/\overline{M}$ is a finitely generated $D/M$-module, that is $[K:k]<\infty$.
\end{proof}

\begin{theorem}\label{006}
(Unequal characteristic case) Let $(D, M)$ be a one-dimensional complete local domain with $\ch D=0$ and $\ch D/M=p>0$. Let $(C, (p))$ be a coefficient ring for $D$. Now $\overline{D}$ is a DVR and $\overline{D}$ is a finitely generated $D$-module, so $0\neq [D:\overline{D}]=\overline{M}^{n}$ for some natural number $n$. Then the following statements are equivalent.
\begin{itemize}
\item[(1)] $D/M\subseteq \overline{D}/\overline{M}$ is a root extension (resp., bounded root extension).
\item[(2)] $C/(p)\subseteq \overline{D}/\overline{M}$ is a root extension (resp., bounded root extension).
\item[(3)] $D\subseteq \overline{D}$ is a root extension (resp., bounded root extension).
\item[(4)] $C+\overline{M}^{n}\subseteq \overline{D}$ is a root extension (resp., bounded root extension).
\item[(5)] $D$ is an AV-domain (resp., API-domain).
\item[(6)] $C+\overline{M}^{n}$ is an AV-domain (resp., API-domain).
\end{itemize}

\end{theorem}
\begin{proof} Now $C/(p)=D/M$, so clearly $(1)\Leftrightarrow (2)$. Also, $(3)\Leftrightarrow (5)$, $(4)\Leftrightarrow (6)$, $(4)\Rightarrow (3)\Rightarrow (1)$. So it is suffice to prove $(2)\Rightarrow (4)$. Let $x\in \overline{D}$, so $(x+\overline{D}/\overline{M})^{k}\in C/(p)$ for some natural number $k$. So there exists $c\in C$ with $x^{k}-c\in\overline{M}$, say $x^{k}=c+m$ where $m\in\overline{M}$. Then $(x^{k})^{p}=(c+m)^{p}=c^{p}+pc^{p-1}m+
(_{2}^{p})c^{p-2}m^{2}+\dots+m^{p}$. Since $p\in \overline{M}$, $x^{kp}-c^{p}\in \overline{M}^{2}$. Continuing we get $x^{kp^{l}}-c^{p^{l}}\in \overline{M}^{l+1}$. Hence $x^{kp^{n-l}}-c^{p^{n-l}}\in \overline{M}^{n}\subseteq D$. Thus $x^{kp^{n-l}}\in C+\overline{M}^{n}$. So $C+\overline{M}^{n}\subseteq \overline{D}$ is a root extension. Also if $C/(p)\subseteq \overline{D}/\overline{M}$ is a bounded root extension, we can take $k$ to be a fixed natural number. Thus $C+\overline{M}^{n}\subseteq \overline{D}$ is a bounded root extension.
\end{proof}

For an integral domain $D$ with quotient field $K$, the \textit{group of divisibility} of $D$ is the Abelian group $G(D):=K^{\ast}/U(D)$ where $U(D)$ is the group of units of $D$, partially ordered by $xU(D)\leq yU(D)$ $\Leftrightarrow$ $Dy\subseteq Dx$. Suppose that $(D,M)$ is one-dimensional local domain. Then $\overline{D}$ is a semilocal PID, so $G(\overline{D})\cong \mathbb{Z}^{n}$ where $n=\mid\max(\overline{D})\mid$. We have a short exact sequence
$$(\ast)\ \ \ \ 0\rightarrow U(\overline{D})/U(D)\rightarrow G(D)\rightarrow G(\overline{D})\rightarrow 0.$$ Since $G(\overline{D})$ is a finitely generated free Abelian group, $(\ast)$ splits. So $G(D)\cong U(\overline{D})/U(D)\oplus G(\overline{D})$. Suppose that $D$ is an AV-domain. Then $\overline{D}$ is a DVR. So $G(\overline{D})\cong \mathbb{Z}$. And since $D\subseteq \overline{D}$ is a root extension, $U(\overline{D})/U(D)$ is torsion. So $G(D)\cong T\oplus \mathbb{Z}$ where $T$ is torsion. If further $D\subseteq \overline{D}$ is a bounded root extension, $U(\overline{D})/U(D)$ is of bounded order, so $G(D)\cong T\oplus \mathbb{Z}$ where $T$ is of bounded order. If $D$ is a DVR, we certainly have $G(D)\cong \mathbb{Z}$. The converse is also true. Suppose that $G(D)\cong \mathbb{Z}$. Then $G(\overline{D})\cong \mathbb{Z}$ and $U(\overline{D})/U(D)=0$. So $\overline{D}$ is a DVR and $U(D)=U(\overline{D})$. Let $x\in \overline{D}$. If $x$ is a unit, then $x\in D$ while if $x$ is not a unit of $\overline{D}$, $1+x$ is a unit of $\overline{D}$, so $1+x\in D$ and hence $x\in D$. So $D=\overline{D}$ is a DVR.
\begin{theorem}\label{007}
Let $(D, M)$ be a one-dimensional local domain with $(\overline{D}, \overline{M})$ a finitely generated $D$-module.
\begin{itemize}
\item[(1)] Suppose that $\ch D/M=0$. Then the following statements are equivalent.
\begin{itemize}
\item[(a)] $D$ is an AV-domain.
\item[(b)] $D$ is an API-domain.
\item[(c)] $D$ is a DVR.
\item[(D)] $G(D)\cong \mathbb{Z}$.
\end{itemize}
\item[(2)]Suppose that $\ch D/M=p>0$. Then the following statements are equivalent.
\begin{itemize}
\item[(a)] $D$ is an AV-domain (resp., API-domain).
\item[(b)] $D\subseteq \overline{D}$ is a root extension (resp., bounded root extension).
\item[(c)] $D/M\subseteq \overline{D}/\overline{M}$ is a root extension (resp., bounded root extension).
\item[(d)] The group of divisibility $G(D)$ of $D$ has the form $G(D)\cong T\oplus \mathbb{Z}$ where $T$ is torsion (resp., of bounded order).
\end{itemize}
\end{itemize}
\end{theorem}
\begin{proof}
Since $\overline{D}$ is a finitely generated $D$-module, by Theorem \ref{004} $D$ is an AV-domain (resp., API-domain) if and only if $\widehat{D}$ is. Also $\widehat{D}/\widehat{M}=D/M$ and $\widehat{\overline{D}}/\widehat{\overline{M}}=\overline{D}/\overline{M}$, so $D/M\subseteq \overline{D}/\overline{M}$ is a (bounded) root extension if and only if $\widehat{D}/\widehat{M}\subseteq \widehat{\overline{D}}/\widehat{\overline{M}}$ is.

$(1)$ \ \ By Theorem \ref{005}, $\widehat{D}$ is an AV-domain if and only if $\widehat{D}$ is a DVR, but $\widehat{D}$ is a DVR if and only if $D$ is a DVR. So $D$ is an AV-domain if and only if $D$ is a DVR and hence $D$ is an API-domain if and only if $D$ is a DVR. The equivalence
of $(c)$ and $(d)$ is given in paragraph preceding Theorem \ref{007}.

$(2)$\ \ $(a)\Leftrightarrow (b)$ The AV-domain case is immediate and does not require $\overline{D}$ to be a finitely generated $D$-module. Now for $\overline{D}$ a finitely generated $D$-module, $[D:\overline{D}]\neq 0$ and hence $D$ is an API-domain if and only if $D\subseteq \overline{D}$ is a bounded root extension (Theorem \ref{04} (4), (5)). Certainly $(b)\Rightarrow (c)$. Suppose that $(c)$ holds, so $\widehat{D}/\widehat{M}\subseteq \widehat{\overline{D}}/\widehat{\overline{M}}$ is a root extension (resp., bounded root extension). First suppose that $\ch D=0$ and $\ch D/M=p>0$. Then Theorem \ref{006} gives that $\widehat{D}$ is an AV-domain (resp., API-domain) and hence so is $D$. Next suppose that $\ch D=\ch D/M$. Then by Theorem \ref{005}, $\widehat{D}$ is an AV-domain (resp., API-domain) and hence so is $D$.

$(b)\Rightarrow(d)$\ \ This follows from the paragraph preceding Theorem \ref{007}.

$(d)\Rightarrow (c)$\ \ Suppose that $G(D)\cong T\oplus \mathbb{Z}$ where $T$ is torsion (resp., of bounded order). Now $G(D)\cong U(\overline{D})/U(D)\oplus\mathbb{Z}^{n}$ where $n=\mid\max(\overline{D})\mid$. Since $T\oplus \mathbb{Z}\cong U(\overline{D})/U(D)\oplus\mathbb{Z}^{n}$, we have $n=1$ and $U(\overline{D})/U(D)$ is torsion (resp., of bounded order). Thus $D/M\subseteq \overline{D}/\overline{M}$ is a root extension (resp., bounded root extension).
\end{proof}

Recall that an integral domain is called a \textit{Cohen-Kaplansky domain} (\textit{CK-domain}) if (1) $D$ is atomic (i.e., each nonzero nonunit is a finite product of irreducible elements) and (2) $D$ has only finitely many irreducible elements up to associates. For results on CK-domains, the reader is referred to \cite{AM}. A domain $D$ is a CK-domain if and only if $D$ is semi(quasi)local and $D$ is locally a CK-domain. A local domain $(D, M)$ is a CK-domain if and only if $D$ is a DVR or $D$ is a one-dimensional analytically irreducible domain with $D/M$ finite. Now if $D$ is a CK-domain, then $D\subseteq \overline{D}$ is a bounded root extension and $\overline{D}$ is a semilocal PID with $|\max{\overline{D}}|=|\max{D}|$. Hence a CK-domain is a API-domain. A local API-domain with $\overline{D}$ a finitely generated $D$-module is a CK-domain if and only if $D$ is a DVR or $D/M$ is finite. Theorem 4.5 \cite{AM} gives a structure theory for complete local CK-domains, compare with Theorem \ref{005} and \ref{006} giving a structure theory for complete local AV-domains and API-domains.

We next consider another type of ``almost DVR". If $D$ is a DVR with a uniformizing parameter $\pi$, then each nonzero nonunit $y$ of $D$ can be written (uniquely) in the form $y=u\pi^{n}$ for some natural number $n$ and unit $u$ of $D$. We next define an ``almost uniformizing parameter".

\begin{definition} \label{05}
Let $D$ be a domain that is not a field.
\rm{An element $a\in D$ is called an \textit{almost uniformizing parameter} for $D$ if for each nonzero nonunit $y\in D$ there exist natural numbers $m$ and $n$ and a unit $u$ of $D$ with $y^{m}=ua^{n}$.}
\rm{A domain having an almost uniformizing parameter is called a \textit{rational almost valuation domain (RAV-domain)}.}
\end{definition}
The name ``rational almost valuation domain" becomes apparent from the next theorem. Note that an almost uniformizing parameter for a domain $D$, if it exists, is far from being unique. Suppose that $a$ is an almost uniformizing parameter for a domain $D$. Let $b$ be a nonzero nonunit of $D$. Then $b^{m}=ua^{n}$ for natural numbers $m$ and $n$ and unit $u$ of $D$. Let $c$ be any nonzero nonunit of $D$. Then $c^{l}=va^{k}$ for natural numbers $l$ and $k$ and unit $v$ of $D$. Then $c^{nl}=v^{n}a^{nk}=v^{n}(u^{-1}b^{m})^{k}=(v^{n}u^{-k})b^{mk}$. So $b$ is also an almost uniformizing parameter for $D$.

\begin{theorem}\label{06}
A domain $D$ is a RAV-domain if and only if $D$ is an AV-domain and $\overline{D}$ is a rational valuation domain (i.e., the value group $G(\overline{D})$ of $\overline {D}$ is order-isomorphic to a subgroup of $(\mathbb{Q}, +)$).
\end{theorem}
\begin{proof}
($\Rightarrow $)\ \ Suppose that $D$ is a RAV-domain. Let $a$ be an almost uniformizing parameter for $D$. We first show that $D$ is an AV-domain. Let $y, z\in D^{\ast}$. If $y$ or $z$ is a unit, then $y\mid z$ or $z\mid y$. So suppose that $y$ and $z$ are nonunits. Then there exist natural numbers $n, n', n, m'$ and units $u, v$ of $D$ with $y^{m}=ua^{n}$ and $z^{m'}=va^{n'}$. Then $y^{mm'}=u^{m'}a^{nm'}$ and $z^{mm'}=v^{m}a^{mn'}$. So $y^{mm'}\mid z^{mm'}$ or $z^{mm'}\mid y^{mm'}$. Hence $D$ is an AV-domain. Thus $\overline{D}$ is a valuation domain and $D\subseteq \overline{D}$ is a root extension. Let $v$ be the valuation associated with $\overline{D}$. Let $z$ be a nonzero nonunit of $\overline{D}$. So there exits a natural number $l$ with $z^{l}\in D$. Since $z^{l}$ is a nonzero nonunit of $D$, there exist natural numbers $m$ and $n$ and a unit $u$ of $D$ with $z^{ml}=ua^{n}$. So $mlv(z)=v(z^{ml})=v(ua^{n})=nv(a)$. Thus $v(z)=\frac{n}{ml}v(a)$ (this makes sense because $G(\overline{D})$ is torsion-free). Let $G=\{q\in\mathbb{Q}\mid qv(a)\in G(\overline{D})\}$. It is easy to check that $G$ is a subgroup of $(\mathbb{Q}, +)$ with $G(\overline {D})=Gv(a)$. It is also easy to check that the map given by $q\rightarrow qv(a)$ is an order-isomorphism. Hence $\overline{D}$ is a rational valuation domain.

$(\Leftarrow)$ Suppose that $D$ is an AV-domain with $\overline {D}$ a rational valuation domain. We may assume that $G(\overline{D})$ is a subgroup of $(\mathbb{Q},+)$. Let $v$ be the valuation associated with $\overline{D}$. Let $a$ be a nonzero nonunit of $D$. We show that $a$ is an almost uniformizing parameter for $D$. Let $b$ be a nonzero nonunit of $D$. Then $v(a)=\frac{m}{n}$ and $v(b)=\frac{m'}{n'}$, where $m,n, m', n'$ are natural numbers. Thus $v(a^{nm'})=mm'=v(b^{mn'})$, so $b^{mn'}=ua^{m'n}$ for some unit $u$ of $\overline{D}$. Since $D\subseteq \overline{D}$ is a root extension, there exists a natural number $k$ with $u^{k}\in D$. So $b^{mn'k}=u^{k}a^{m'nk}$ where $u^{k}$ is a unit of $D$. Hence $a$ is an almost uniformizing parameter for $D$, so $D$ is a RAV-domain.
\end{proof}

The title of this paper is Almost Discrete Valuation  Domains. So what is (or what should be) an almost DVR? We have investigated several generalizations: (quasi)local AV-domains with $\overline{D}$ a DVR, (quasi)local API-domains, and RAV-domains. Now we consider the following seven conditions for a (necessarily) quasilocal integral domain $(D,M)$:
\begin{itemize}
\item[(1)] $D$ is a local API-domain,
\item[(2)] $D$ is a Noetherian AV-domain,
\item[(3)] $D$ is an AV-domain and for each ideal $A$ with $0\subsetneq A\subsetneq D$, there exists a natural number $n$ with $M^{n}\subseteq A$,
\item[(4)] $D$ is an AV-domain with $\overline{D}$ a DVR,
\item[(5)] $D$ is an AV-domain with $\bigcap_{n=1}^{\infty}M^{n}=0$,
\item[(6)] $D$ is a quasilocal API-domain, and
\item[(7)]  $D$ is a RAV-domain.
\end{itemize}

Our next theorem gives some of the implications among $(1)$-$(7)$.
\begin{theorem} \label{07}
Let $(D,M)$ be a quasilocal integral domain.
\begin{itemize}
\item[(a)] We have the following implications

$(1)\Rightarrow (2)\Rightarrow (3)\Rightarrow (5)$,
$(2)\Rightarrow (4)\Rightarrow (5)$,
$(1)\Rightarrow (6)$, and
 $(4)\Rightarrow (7)$.
\item[(b)] We also have

$(2)\nRightarrow (6)$ and $(6)\nRightarrow (2)$ (so $(2)\nRightarrow (1)$ and $(6)\nRightarrow(1)$),

$(3)\nRightarrow(2)$,  $(4)\nRightarrow(2)$ and

$(7)\nRightarrow (5)$ (so $(7)\nRightarrow (3)$ and $(7)\nRightarrow (4)$).
\end{itemize}
\end{theorem}
\begin{proof}
$(a)$\ \ $(1)\Rightarrow (2)$\ \ Clear.

$(2)\Rightarrow (3), (4)$\ \ Since $D$ is a Noetherian AV-domain, $\overline{D}$ is a Krull domain. Since $\overline{D}$ is a valuation domain, $\overline{D}$ is a DVR. This gives (4). Now $\overline{D}$ a DVR gives $\dim{D}=1$. So $D$ is a one-dimensional local domain and hence (3) holds.

$(3)\Rightarrow (5)$\ \ Here $\bigcap_{n=1}^{\infty}M^{n}\subseteq\bigcap\{A\mid 0\subsetneq A\subsetneq D$ is an ideal of $D\}=0$.

$(4)\Rightarrow (5)$\ \ Let $Q$ be the maximal ideal of $\overline{D}$. So
$\bigcap_{n=1}^{\infty}M^{n}\subseteq\bigcap_{n=1}^{\infty}Q^{n}=0$.

$(1)\Rightarrow (6)$\ \ Trivial.

$(4)\Rightarrow (7)$\ \ Suppose that $\overline{D}$ is a DVR. Then $\overline{D}$ is a rational valuation domain, so $D$ is a RAV-domain by Theorem \ref{06}.

$(b)$\ \ This will follow from  Example \ref{09}.
\end{proof}

The following diagram summarizes Theorem \ref{07}.

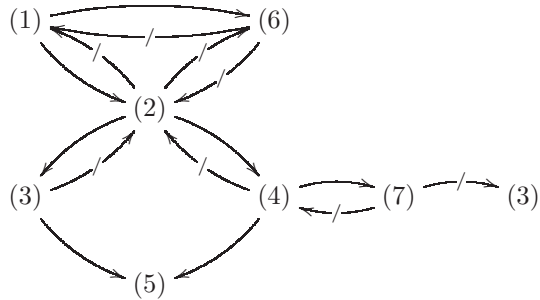
$$\xymatrix@R=17pt@C=28pt{
 & (1)\ar@/_/[rd]|{}\ar@/^/[rr]|{} & &(6)\ar@/^/[ll]|{/}\ar@/^/[dl]|{/}\\
& & (2)\ar@/_/[lu]|{/}\ar@/^/[ru]|{/}\ar@/_/[dl]\ar@/^/[dr]\\
&(3) \ar@/_/[ur]|{/}\ar@/_/[dr]& &(4)\ar@/^/[ul]|{/}\ar@/^/[dl] \ar@/^/[r]&(7) \ar@/^/[l]|{/}\ar@/^/[r]|{/}&(3)\\
& &(5) }$$

Our next theorem summarizes when the $D+M$ construction satisfies one of the conditions $(1)$-$(7)$. This will be used to show that certain implications given in Theorem \ref{07} cannot be reversed.
\begin{theorem}\label{08}
Let $(V,M)$ be a valuation domain with the quotient field $L$ of the form $V=K+M$ (e.g., $V=K[X]_{(f)}$ where $f$ is linear or $V=K[[X]]$, $K$ is a field). Let $D$ be a subring of $K$ having quotient filed $F$, $\overline{D}$ (resp., $\widetilde{D}$) the integral closure of $D$ in $F$ (resp., $K$) and $R=D+M$. So $R$ is a subring of $V$ with quotient field $L$.
\begin{itemize}
\item[(1)] $R=V$ $\Leftrightarrow$ $D=K$. If $R\subsetneq V$, $[R:V]=M$.
\item[(2)] $\overline{R}=\widetilde{D}+M$.
\item[(3)] $R$ is quasilocal $\Leftrightarrow$ $D$ is quasilocal.
\item[(4)] $R$ is local $\Leftrightarrow$ $R$ is Noetherian $\Leftrightarrow$ $V$ is a DVR, $D=F$ is a field and $[K:F]<\infty$.
\item[(5)] $R$ is a valuation domain $\Leftrightarrow$ $D$ is a valuation domain and $K=F$.
\item[(6)] $R$ is a DVR (resp., rational valuation domain) $\Leftrightarrow$ $D=K$ and $V$ is a DVR (resp., rational valuation domain) $\Leftrightarrow$ $R=V$ is a DVR (resp., rational valuation domain).
\item[(7)] $R\subseteq \overline{R}$ is a (bounded) root extension $\Leftrightarrow$ $D\subseteq \widetilde{D}$ is a (bounded) root extension.
\item[(8)] $R$ is an AV-domain $\Leftrightarrow$
$\widetilde{D}$ is a valuation domain with quotient field $K$ and $D\subseteq \widetilde{D}$ is a root extension
$\Leftrightarrow$ $D$ is an AV-domain and $F\subseteq K$ is a root extension (i.e., $F\subseteq K$ is a purely inseparable extension or $K$ is algebraic over a finite field.).

\item[(9)] $R$ is an AV domain with $\overline{R}$ a DVR (necessarily $\overline{R}=V$) $\Leftrightarrow$ $(R, Q)$ is an AV-domain with $\bigcap_{n=1}^{\infty}Q^{n}=0$ $\Leftrightarrow$ $(R, Q)$ is an AV-domain and for each ideal $A$ of $R$ with $0\subsetneq A\subsetneq R$ there exists a natural number $n$ with $Q^{n}\subseteq A$ $\Leftrightarrow$
 $V$ is a DVR, $D=F$ is a field and $F\subseteq K$ is a root extension.

\item[(10)] $R$ is a local AV-domain $\Leftrightarrow$ $V$ is a DVR, $D=F$ is a field with $[F:K]<\infty$ and $F\subseteq K$ is a root extension.
\item[(11)] $R$ is a quasilocal API-domain $\Leftrightarrow$ $V$ is a DVR, $D=F$ is a field and $F\subseteq K$ is a bounded root extension (i.e., $F\subseteq K$ is purely inseparable of bounded exponent or $K$ is finite).
\item[(12)] $R$ is a local API-domain $\Leftrightarrow$ $V$ is a DVR, $D=F$ is a field with $[F:K]<\infty$, and $F\subseteq K$ is a bounded root extension.
\item[(13)] $R$ is a RAV-domain $\Leftrightarrow$ $V$ is a rational valuation domain, $D=F$ is a field, and $F\subseteq K$ is a root extension.
\end{itemize}
\end{theorem}
\begin{proof}
$(1)$-$(5)$ are well-known.

$(6)$ \ \ This follows from (5) and the facts that $\dim R=\dim D+\dim V$ and a rational valuation domain (and hence a DVR) is one-dimensional.

$(7)$\ \ Clear.

$(8)$\ \ By Theorem \ref{02} (1), $R$ is an AV-domain if and only if $\overline{R}=\widetilde{D}+M$ is a valuation domain (equivalently, $\widetilde{D}$ is a valuation domain with the quotient field $K$) and $R\subseteq \overline {R}$ is a root extension (equivalently, $D\subseteq \widetilde{D}$ is a root extension). Suppose that $\widetilde{D}$ is a valuation domain with quotient field $K$ and $D\subseteq \widetilde{D}$ is a root extension. By Theorem \ref{02} (4), $D$ is an AV-domain. And since  $D\subseteq \widetilde{D}$ is a root extension, $L=T(D)\subseteq T(\widetilde{D})=K$ is a root extension. Conversely, suppose that $D$ is an AV-domain and $F\subseteq K$ is a root extension. Then $D\subseteq \widetilde{D}$ is a root extension. For if $x\in \widetilde{D}\subseteq K$, then $x^{n}\in F$ for some natural number $n$. But $x^{n}$ is integral over $D$, so $x^{n}\in \overline{D}$. But $D\subseteq \overline{D}$ is also a root extension, so there exists a natural number $m$ with $x^{mn}=(x^{n})^{m}\in D$. So $D\subseteq \widetilde{D}$ is a root extension. The condition equivalent to $F\subseteq K$ being a root extension is Theorem \ref{01} (1).

$(9)$\ \ $R$ is an AV-domain with $\overline{R}$ a DVR (equivalently, $\widetilde{D}=K$ and $\overline{R}=V$ is a DVR) $\Leftrightarrow$ $D=F$ is a field with $F\subseteq K$ a root extension and $V$ is a DVR. Suppose that the last condition holds. Let $M=V\pi$ be the maximal ideal of $V$. So $M$ is also the maximal ideal of $R$. Let $A$ be a nonzero proper ideal of $R$. Now some $\pi^{n}\in R$ since $R\subseteq \overline{R}$ is a root extension and hence $\pi^{mn}=(\pi^{n})^{m}\in A$ since $\dim R=1$, and so $\sqrt{A}=M$. Then $M^{mn+1}=\pi^{mn}M\subseteq AM\subseteq A$. Finally, suppose that $(R,Q)$ is an AV-domain with $\bigcap_{n=1}^{\infty}Q^{n}=0$. Now $M\subseteq Q$, so $\bigcap_{n=1}^{\infty}M^{n}\subseteq\bigcap_{n=1}^{\infty}Q^{n}=0$ and hence $V$ is a DVR. Now if $D$ is not a field, then $Q\supsetneq M$ and hence $Q^{n}\supseteq M$ for each natural number $n$. Then $M\subseteq \bigcap_{n=1}^{\infty}Q^{n}=0$, a contradiction.

$(10)$\ \ This follows from (4) and (8).

$(11)$\ \ ($\Rightarrow$) Suppose that $R$ is a quasilocal API-domain. Since $[R:V]=M\neq 0$, $R\subseteq \overline{R}$ is a bounded root extension and $\overline{R}$ is a DVR by Theorem \ref{04} (5). By (9), $V$ is a DVR and $D=F$ is a field. Since $R\subseteq \overline{R}$ is a bounded root extension, $F\subseteq K$ is a bounded root extension.

($\Leftarrow$) $\overline{R}=V$ is a DVR and $R\subseteq \overline{R}$ is a bounded root extension since $D=F\subseteq K$ is a bounded root extension. By Theorem \ref{04} (4), $R$ is a a quasilocal API-domain.

$(12)$\ \ This follows from (4) and (11).

$(13)$\ \ ($\Rightarrow$) Suppose that $R$ is a RAV-domain. Then $\overline{R}$ is a rational valuation domain by Theorem \ref{06}. By (6), $\overline{R}=V$, so $V$ is a rational valuation domain. Since $\dim{R}=\dim{\overline{R}}=\dim{V}=1$, $D=F$ is a field. By (8), $F\subseteq K$ is a root extension.

($\Leftarrow$) Here $\overline{R}=V$ is a rational valuation domain, so by Theorem \ref{06}, $R$ is a RAV-domain.
\end{proof}

\begin{example}\label{09}
(Non-implications of Theorem \ref{07})
\begin{itemize}
\item[(a)] $(2)\nRightarrow (6)$ (local AV-domain $\nRightarrow$ quasilocal API-domain)
 \rm{ \cite[Example 3.6]{AKL} Let $F=\bigcup_{n=1}^{\infty} GF(p^{2^{n}})$ and $K=F(GF(p^{3}))$. So $[F:K]<\infty$ and $F\subsetneq K$ is a root extension but not a bounded root extension (Theorem \ref{01}). So $R=F+XK[[X]]$ is a local AV-domain but not an API-domain by Theorem \ref{08} (12) or Theorem \ref{005}.}

\item[(b)] \rm{$(6)\nRightarrow (2)$ (quasilocal API-domain $\nRightarrow$ local AV-domain) \cite[Example 3.7]{AKL} Let $\{X_{n}\}_{n=1}^{\infty}$ be a set of indeterminates over $\mathbb{Z}_{p}$. Let $K=\mathbb{Z}_{p}(\{X_{n}\}_{n=1}^{\infty})$ and $F=\mathbb{Z}_{p}(\{X_{n}^{p}\}_{n=1}^{\infty})$. So $F\subset K$ is a infinite dimensional bounded root extension. Hence $R=F+XK[[X]]$ is a non-Noetherian quasilocal API-domain by Theorem \ref{08} (10), (11).}

\item[(c)] $(3)\nRightarrow (2), (4)\nRightarrow (2)$ ($(R,Q)$ an AV-domain with $Q^{n}\subseteq A$ for each ideal $0\subsetneq A\subsetneq R$, or $\overline{R}$ a DVR $\nRightarrow$ $R$ is a local AV-domain)\rm{ Use (a) and apply Theorem \ref{08} (9), (10).}

\item[(d)]$(7)\nRightarrow (5)$ ($R$ is a RAV-domain $\nRightarrow$ $(R,Q)$ is an AV-domain with $\bigcap_{n=1}^{\infty}Q^{n}=0$) \rm{Just take a valuation domain with value group $(\mathbb{Q},+)$.}
\end{itemize}
\end{example}

\end{document}